\documentclass[12pt]{article}
\usepackage{sw20elba}
\usepackage{amssymb,amsmath}

\newtheorem{theorem}{Theorem}[section]

\newtheorem{definition}[theorem]{Definition}

\newtheorem{proposition}[theorem]{Proposition}
\newtheorem{remark}[theorem]{Remark}

\newenvironment{proof}[1][Proof]{\textbf{#1.} }{\ \rule{0.5em}{0.5em}}

\begin{document}

\author{Robert Erdahl and Konstantin Rybnikov}
\title{An Infinite Series of Perfect Quadratic Forms and Big Delaunay Simplexes in $%
\mathbb{Z}^{n}$.}
\date{May 17, 2001}
\maketitle

\begin{abstract}
George Voronoi (1908-09) introduced two important reduction methods for
positive quadratic forms: the reduction with perfect forms, and the
reduction with $L$-type domains.

A form is \emph{perfect }if can be reconstructed from all representations of
its arithmetic minimum. Two forms have the same $L$-type if Delaunay tilings
of their lattices are affinely equivalent. Delaunay (1937-38) asked about
possible relative volumes of lattice Delaunay simplexes. We construct an
infinite series of Delaunay simplexes of relative volume $n-3$, the best
known as of now. This series gives rise to a new infintie series of perfect
forms $TF_{n}$ with remarkable properties: e.g. $TF_{5}=D_{5},$ $%
TF_{6}=E_{6}^{\ast },$ $TF_{7}=\varphi _{15}^{7};$ for all \ $n$ the domain
of $TF_{n}$ is adjacent to the domain of the $2$-nd perfect form $D_{n}.$
Perfect form $TF_{n}$\ is a direct $n$-dimensional generalization of Korkine
and Zolotareff's $3$-rd perfect form $\phi {2}^{5}$ in $5$ variables. It is
likely that this form is equivalent to Anzin's (1991) form $h_n$.
\end{abstract}

\smallskip


 \smallskip

\section{Introduction and main result\label{intro}}

Positive quadratic forms (referred to as PQFs) in $n$ indeterminate form a 
\emph{closed }cone $\frak{P}(n)$ of dimension $N=\frac{n(n+1)}{2}$ in $%
\mathbb{R}^{N}$, and this cone is the main object of study in our paper. The
interior of $\frak{P}(n)$ consists of positive definite forms of rank $n$.
We abbreviate positive \emph{definite} forms as PDQFs. PDQFs serve as
algebraic representations of \textit{point lattices}. There is a one-to-one
correspondence between isometry classes of $n$-lattices and integral
equivalence classes (i.e. with respect to $GL(n,\mathbb{Z})$-conjugation) of
PDQFs in $n$ indeterminates. For basic results of the theory of lattices and
PQFs and their applications see \textit{Ryshkov and Baranovskii }(1978), 
\textit{Gruber and Lekkerkerker} (1987), \textit{Erd\"{o}s, Gruber and Hammer%
} (1989), \textit{Conway and Sloane }(1999).

$GL(n,\mathbb{Z})$ acts pointwise on the space of quadratic forms $Sym(n,%
\mathbb{R)}\cong \mathbb{R}^{N}$. A \emph{polyhedral reduction partition} $R$
of $\frak{P}(n)$ is a partition of this cone into open \emph{convex
polyhedral} cones such that:

\begin{definition}
\begin{enumerate}
\item  it is invariant with respect to $GL(n,\mathbb{Z})$;

\item  there are finitely many inequivalent cones in this partition;

\item  for each cone $C$ of $R$ and any PQF $f$ in $n$ indeterminates, $f$
can be $GL(n,\mathbb{Z})$-equivalent to at most finitely many forms lying in 
$C$.
\end{enumerate}
\end{definition}

The \emph{partition into perfect cones} and the \emph{$L$-type partition}
(also referred to as the Voronoi partition of the 2nd kind, or the partition
into Voronoi reduction domains) are important polyhedral reduction
partitions of $\frak{P}(n)$. (Our usage of term \emph{domain }is lax; it
should be clear from the context whether we mean the whole arithmetic class,
or just \ one element of this class.) These partitions have been intensively
studied in geometry of numbers since times of \textit{Korkin, Zolotareff}
(1873) and \textit{Voronoi} (1908-1909), and more recently in combinatorics
(e.g. \textit{Deza} \textit{et} \textit{al}. (1997)), and algebraic geometry
(e.g. \textit{Alexeev} (1999a,b)). \emph{In most previous works} (e.g. 
\textit{Voronoi }(1908,1909), \textit{Ryshkov, Baranovskii} (1976), \textit{%
Dickson} (1972)) \emph{the }$L$\emph{-type partition of }$P(n)$\emph{, or
sub-cones of }$P(n)$\emph{, was constructed by refining the perfect
partition.} It is not an exaggeration to say that in almost any systematic
study, except for Engels' computational investigations, $L$-types were
approached via perfect forms. For example, \textit{Voronoi }started
classifying 4-dimensional $L$-types by analyzing the Delaunay ($L$-)tilings
of forms lying in the 1st ($A_{n}$) and 2nd ($D_{n}$) perfect domains. The
same route was followed by \textit{Ryshkov and Baranovskii }(1976). It was
wideley believed that the $L$-type parition is the refinement of the perfect
parition, i.e. each convex cone of the perfect partition is the union of
finitely many conex cones from the $L$-type partition. This conjecture is
implicit in \textit{Voronoi}'s memoirs (1908-1909), and explicit in \textit{%
Dickson }(1972), where he showed that the first perfect domain is the only
perfect domain which coincided with an $L$-type domain. \ More information
on this refinement conjecture and its failure for $\frak{P}(n)\subset Sym(6,%
\mathbb{R)}$ can be found in our note ''Voronoi-Dickson hypothesis on
perfect forms ands L-types'' published in this volume.

In our paper we continue the above-mentioned line of research on the
relashionship between these two reduction partitions. However, we go in the 
\emph{opposite} direction. We construct an arithmetic class of perfect forms
from an arithmetic class of $L$-types.

\begin{remark}
\textbf{Since there is a lot of overlap in references between this paper and our paper on 
 the Voronoi-Dickson hypothesis, some of the references for this paper should be found in the
bibliography for the other paper. Such references are
quoted in italics, e.g. } Baranovski (1991).
\end{remark}

\section{Perfect and $L$-type partitions}

\subsection{$L$-types}

\begin{definition}
Let $L$ be a lattice in $\mathbb{R}^{n}$. A convex polyhedron $P$ in $%
\mathbb{R}^{n}$ is called a Delaunay cell of $L$ with respect to a positive
quadratic form $f\mathbf{(x,x)}$ if:

\begin{enumerate}
\item  for each face $F$ of $P$ we have $conv(L\cap F)=F$;

\item  there is a quadric circumscribed about $P$, called the \emph{empty
ellipsoid} of $P$ when $f\mathbf{(x,x)}$ is positive definite, whose
quadratic form is $f\mathbf{(x,x)}$ (in case $rank\ f<n$, this quadric is an
elliptic cilinder);

\item  no points of $L$ lie inside the quadric circumscribed \ about $P$.
\end{enumerate}
\end{definition}

When $f=\sum_{i=1}^{n}x_{i}^{2}$, our definition coincides with the
classical definition of Delaunay cell in $\mathbb{E}^{n}$. Delaunay cells
form a convex face-to-face tiling of $L$ that is uniquely defined by $L$ (%
\textit{Delaunay}, 1937). Two Delaunay cells are called \emph{homologous} is
they can be mapped to each other with a composition of a lattice translation
and a central inversion with respect to a lattice point.

\begin{definition}
PQFs $f_{1}$ and $f_{2}$ belong to the same convex $L$-domain if the
Delaunay tilings of $\mathbb{Z}^{n}$ with respect to $f_{1}$ and $f_{2}$ are
identical. $f_{1}$ and $f_{2}$ belong to the same $L$-type if these tilings
are equivalent with respect to $GL(n,\mathbb{Z)}$.
\end{definition}

The folowing proposition establishes the equivalence between the Delaunay's
definition of $L$-equivalence for lattices and the notion of $L$-equivalence
for arbitrary PQFs, which is introduced above.

\begin{proposition}
Positive definite forms $f$ and $g$ \ belong to the same $L$-type if the
corresponding lattices belong to the same\ $L$-type with respect to the form 
$\sum_{i=1}^{n}x_{i}^{2}$.
\end{proposition}

\begin{theorem}
(Voronoi) The parition of $\frak{P}(n)$ into $L$-types is a reduction
partition. Moreover, it is face-to-face.
\end{theorem}

The notions of Delauny tiling and $L$-type are important in the study of
extremal and group-theoretic properties of lattices (see also
''Voronoi-Dickson Hypothesis...'').  Barnes and Dickson (1967, 1968) and,
later, in a geometric form, \textit{Delaunay }et al. (1969, 1970) proved the
following

\begin{theorem}
The \emph{closure} of any $N$-dimensional convex $L$-domain contains at most
one local minimum of the sphere covering density. The group of $GL(n,\mathbb{%
Z)}$-automorphisms of the domain maps this form to itself.
\end{theorem}

Using this approach, \textit{Delaunay, Ryshkov and Baranovskii }(1963, 1976)
found the best lattice coverings in $\mathbb{E}^{4}$ and $\mathbb{E}^{5}.$
The theory of $L$-types also has numerous connections to combinatorics and,
in particular, to cuts, hypermetrics, and regular graphs (see \textit{Deza }%
et al. (1997)). Recently, \textit{V. Alexeev} (1999a,b) found exciting
connections between compactifications of moduli spaces of principally
polirized abelian varities and $L$-types.

\subsection{Perfect cones}

The $L$-type partition of $\frak{P}(n)$ is closely related to the theory of 
\emph{perfect forms} originated by \textit{Korkine and Zolotareff} (1873).
Let $f\mathbf{(x,x})$ be a PDQF. The \emph{arithmetic minimum } of $f\mathbf{%
(x,x})$ is the minimum of this form on $\mathbb{Z}^{n}$. The integral
vectors on which this minimum is attained are called the representations of
the minimum, or the \emph{minimal vectors} of $f\mathbf{(x,x})$: these
vectors have the minimal length among all vectors of $\mathbb{Z}^{n}$ when $f%
\mathbf{(x,x})$ is used as the metrical form. Form $f\mathbf{(x,x})$ is
called \emph{perfect} if it can be reconstructed up to scale from all
representations of its arithmetic minimum. In other words, a form $f\mathbf{%
(x,x})$ with the arithmetic minimum $m$ and the set of minimal vectors $\{%
\mathbf{v}_{k}|\,k=1,...,2s\}$ is perfect if the system 
\begin{tabular}{|l|}
\hline
$\sum_{i,j=1}^{n}a_{ij}v_{k}^{i}v_{k}^{j}=m,$ \\ \hline
\end{tabular}
where $k=1,...,2s,$ has a unique solution $(a_{ij})$ in $Sym(n,\mathbb{R)}%
\cong \mathbb{R}^{N}$ (indeed, uniquness requires at least $n(n+1)$ minimal
vectors).

\begin{definition}
PQFs $f_{1}$ and $f_{2}$ belong to the same cone of the perfect partition if
they both can be written as \emph{strictly positive} linear combinations of
some subset of minimal vectors of a perfect form $\phi $. $f_{1}$ and $f_{2}$
belong to the same perfect type if there is $f_{1}^{\prime }$, equivalent to 
$f_{1}$, such that $f_{1}^{\prime }$ and $f_{2}$ belong to the same cone of
the perfect partition.
\end{definition}

\begin{theorem}
(Voronoi) The partition of $\frak{P}(n)$ into perfect domains is a reduction
partition. Moreover, it is face-to-face. Each 1-dimensional cone of this
partition lies on $\partial \frak{P}(n)$.
\end{theorem}

Perfect forms play an important role in lattice sphere packings. Voronoi's
theorem (1908) says that if a form is\emph{\ extreme}---i.e., a maximum of
the packing density---it must also be perfect (see \textit{Coxeter }(1951), 
\textit{Conway, Sloane} (1988) for the proof). The notion of \emph{eutactic
form} arises in the study of the dense lattice sphere packings and is
directly related to the notion of perfect form. The reciprocal\ of $f(%
\mathbf{x,x})$ is a form whose Gramm matrix is the inverse of the Gramm
matrix of $f(\mathbf{x,x}).$ The dual form is normally denoted by $f^{\ast }(%
\mathbf{x,x}).$ A form $f(\mathbf{x,x})$ is called \textit{eutactic }if the
dual form $f^{\ast }(\mathbf{x,x})$ can be written as $\sum_{k=1}^{s}\alpha
_{k}(\mathbf{v}_{k}\bullet \mathbf{x)}^{2},$ where $\{\mathbf{v}%
_{k}|\,k=1,...,s\}$ is the set of mutually non-collinear minimal vectors of $%
f(\mathbf{x,x}),$ and $\alpha _{k}>0$.

\begin{theorem}
(Voronoi) A\ form $f(\mathbf{x,x})$ is a maximum of the sphere packing
density if and only if $f(\mathbf{x,x})$ is perfect and eutactic.
\end{theorem}

Voronoi gave an algorithm finding all perfect domains for given $n$. This
algorithm is known as Voronoi's reduction with perfect forms. \ For the
computational analysis of his algorithm and its improvements see \textit{%
Martinet }(1996). The perfect forms and the incidence graphs of perfect
paritions of $\frak{P}(n)$ have have been completely described for $n\leq 7$

\subsection{Relashionship between perfect domains and L-types: interpretable
and non-interpretable perfect walls.}

\textit{Voronoi }(1908-09) proved that for $n=2,3$ the $L$-partition and the
perfect partition of $\frak{P}(n)$ coincide. The perfect facet $D_{4}$ (the
2nd perfect form in 4 variables) exemplifies a new pattern in the relation
of these partitions. Namely, the facet $D_{4}$ is decomposed into a number
of simplicial $L$-type domains like a pie:\ this decomposition consists of
the cones with apex at the affine center of this facet over the $(N-2)$%
-faces. These simplexes are $L$-type domains of two arithmetic types: type I
is adjacent to the the perfect/$L$-type domain of $A_{4},$ type II is
adacent to an arithmetically equivalent $L$-type domain (also type II,
indeed) from the $L$-subdivision of the adjacent $D_{4}$ domain (for details
see \textit{Delaunay et al.} (1963, 1968)).

Following the lead of Voronoi, \textit{Delaunay et al.} proved that for $n=4$
the tiling of $\frak{P}(n)$ with $L$-type domains refines the partition of
this cone into perfect domains. \textit{Ryshkov and Baranovskii }(1975)
proved the refinement hypothesis for $n=5$. See another our paper from this
volume to learn why this hypothesis fails for $n=6$. In cases where L-type
domains refine perfect ones, the L-type is changing on each perfect wall. We
call such perfect walls\emph{\ interpretable}. It is not yet clear why some
perfect walls are intepretable, while other, like the wall between domains
of types $E_{6}$ and $E_{6}^{\ast },$ are not.

Below we construct an infinite series of Delaunay polytopes $R_{n}$ on $n+2$
vertices in $\mathbb{Z}^{n}$ (Theorem \ref{Construction of Rep Complex}$.$
We prove that they are repartitoning complexes. One of the two
triangulations of this polytope has a Delaunay simplex of relative volume $%
n-3$. This triangulation defines an $N$-dimensional $L$-type domain which is
a subcone of a perfect domain $TF_{n},$\ that we describe in Theorem \ref
{Perfect Domain Description}. Domain of type $TF_{n}$ shares a wall with
domain of type $D_{n}.$ All forms lying on this wall have the repartitiong
complex $R$ in its Delaunay tiling.

\section{Fat Symplexes}

It is well known that the Delaunay tiling of lattice $E_{6}$ consists of
Gosset polytopes $G$ (e.g. see \textit{Baranovskii} 1991). For properties of 
$G$-tope see \textit{Coxeter }(1973,1995). All faces of $G$-tope are regular
polytopes. Let $S_{4}$ be a 4-face of $G$-tope which is a common facet of
two cross-polytopal facets. Since any pair of vertices of $G$-tope is either
a diagonal or an edge of a cross-polytopal facet, there are only two
vertices of the polytope which do not have common edges with vertices of $%
S_{4}$. The volume of the convex hull of $S_{4}$ and these two ``distant''
vertices is 3 times the volume of a fundamental simplex of $E_{6}$ . Using a
computer program we checked that there are no other simplexes of volume 3,
Delaunay or not, in $G$-tope (below we prove that all simplexes inscribed
into $G$-tope are Delaunay). There are exactly 216 4-faces that serve as
common facets of pairs of cross-polytopal facets and all of them are
equivivalent with respect to the group of the $G$-tope. Therefore, there are
exactly 216 Delaunay simplexes of relative volume 3 (we will often omit the
word \emph{relative}) in $G$-tope. They are all equivalent with respect to
the isometry group of $G$-tope. According to Ryshkov and Baranovskii (1998)
there is only one arithmetic type of triple Delaunay simplexes in $6$%
-lattices, and there are no Delaunay simplexes of volume greater than $3$ in 
$6$-dimensional lattices.

In an appropriate coordinate system the vertices (here the column-vectors)
of this simplex have the following form.

\begin{equation*}
\left( 
\begin{array}{ccccccc}
0 & 1 & 0 & 0 & 0 & 0 & 1 \\ 
0 & 0 & 1 & 0 & 0 & 0 & 1 \\ 
0 & 0 & 0 & 1 & 0 & 0 & 1 \\ 
0 & 0 & 0 & 0 & 1 & 0 & 1 \\ 
0 & 0 & 0 & 0 & 0 & 1 & 1 \\ 
0 & 0 & 0 & 0 & 0 & 0 & -3
\end{array}
\right)
\end{equation*}

There is a repartitioning complex $R_{6}$ one of whose triangulations
includes the above simplex. The vertices of this repartitioning complex are
the vertices of the above simplex plus the vertex $(0,0,0,0,0,0,1)^{T}$.

\textit{Delaunay }(1937) asked about possible volumes of Delaunay simplexes.
Ryshkov (1973) showed that in every dimension $2r+1$ there is a lattice with
a Delaunay simplex of relative volume $r$. Namely, Ryshkov proved that
lattice $A_{n}^{k}$ for \ \ \ \ $n\geq 2k+1$ has a Delaunay simplex of rel.
volume $k$. Ryshkov also noticed that in the case of $A_{n}^{k}$ the
existence of big Delaunay simplexes is closely related to another
interesting phenomenon: for $n\geq $ $9$ perfect lattice $A_{n}^{k}$ is not
generated by its shortest vectors Earlier, \textit{Coxeter} (1951) made a
similar observation about the relevance of these two phenomena in case of $%
A_{n}^{k},$ but he did not know for sure if $A_{n}^{k}$ had such big
simplexes.

We generalized the construction of the above simplex to the following series
of simplexes of volume $n-3$. Although, to our knowledge, this is the best
infinite series of big Delaunay simplexes, in Leech lattice $\Lambda _{24}$
all Delaunay simplexes are non-fundamental, and the biggest of them has
volume 20480. Haase and Ziegler (2000) showed that for $n>3$ there are empty
lattice simplexes of arbitrary large volume (not Delaunay, indeed). A\
trivial upper bound on the rel. volume of a Delaunay simplex is $\dfrac{n!}{2%
}.$

\begin{equation*}
S_{n}=\left( 
\begin{array}{ccccccc}
0 & 1 & 0 & 0 & 0 & \dots & 1 \\ 
0 & 0 & 1 & 0 & 0 & \dots & 1 \\ 
0 & 0 & 0 & 1 & 0 & \dots & 1 \\ 
0 & 0 & 0 & 0 & 1 & \dots & 1 \\ 
\dots & \dots & \dots & \dots & \dots & \dots & \dots \\ 
0 & 0 & 0 & 0 & 0 & \dots & -(n-3)
\end{array}
\right)
\end{equation*}
The corresponding repartitioning complex is obtained by adjoining vertex $%
(0,...,0,1)^{T}$. In this paper we use a short-hand notation for $n$-vectors
that have few distinct integral coordinates and for families of such vectors
obtained from some $n$-vector by all circular permutations of selected
subsets of its coordinates. Here are the rules:

\begin{enumerate}
\item  $m^{k}$ stands for $k$ consecutive positions filled with $m$'s.

\item  Square brackets $[a_{1}...a_{n}]$ are used to denote a vector (or
type of vector) that can be obtained from this vectors by circular
permutations in sequences of coordinates that are separated by commas and
bordered on the sides by semicolons and/or brackets.

\item  A family of vectors that are obtained from vector $%
[a_{1}...a_{n}]=(a_{1}...a_{n})^{T}$ by all admissible (see 2) permutations
is denoted by $[a_{1}...a_{n}]^{\#},$ where $\#$ is the number of such
vectors; if $\#=1$, we omit $\#.$
\end{enumerate}

Example: $[1^{n-3},0^{2};3]^{\binom{n-1}{2}}$ stands for all vectors with $3$
at the last entry, two $0$'s and $n-3$ $1$'s among the first $n-1$
coordinates.

To prove that the above simplex is a lattice Delaunay simplex we need the
theory of $(0,1)$-dual systems developed by Erdahl and Ryshkov (1990, 1991
a,b). Let $S$ be a set of integer vectors in $\mathbb{Z}^{n}$. The $(0,1)$%
-dual of $S$ is the set of all integer vectors in $\mathbb{Z}^{n}$ that have
the scalar product of 0 or 1 with all vectors of $S$. We denote the $(0,1)$%
-dual of $S$ by $S^{0}$. Erdahl and Ryshkov (1990) showed that if the double
dual of an integral simplex has only $n+2$ points, then this simplex is a
Delaunay simplex for some PDQFs.

\begin{theorem}
(Erdahl, Ryshkov) Let $S$ be a set of vectors in $\mathbb{Z}^{n}$. If $%
(S^{0})^{0}\backslash $ $[0^{n}\mathbf{]}$ consists of $n+1$ linearly
independent vectors, then there is an N-dimensional cone of PDQFs for which $%
S$ is a Delaunay simplex in $\mathbb{Z}^{n}$.
\end{theorem}

Using the Erdahl-Ryshkov theorem we verify in the following proposition that
our series is indeed a series of lattice Delaunay simplexes.

\begin{theorem}
\label{Construction of Rep Complex}For any $n>3$ there are lattices with
Delaunay simplex of volume $n-3$.\emph{\ }
\end{theorem}

\begin{proof}[Proof]
Let $S_{n}$ be a simplex in $\mathbb{Z}^{n}$ whose vertices are the columns
of the following matrix: 
\begin{equation*}
S_{n}=\left( 
\begin{array}{cc}
I_{n-1} & \mathbf{1}_{n-1} \\ 
\mathbf{0}_{n-1}^{T} & -(n-3)
\end{array}
\right) .
\end{equation*}
Here $I_{n-1}$ is the $(n-1)\times (n-1)$ identity matrix, $1_{n-1}$ is the
column $(n-1)$-vector of ones, and $0_{n-1}^{T}$ is the row $(n-1)$-vector
of zeros. $S_{n}^{0}$ cannot have vectors with negative numbers in positions 
$1$\ through $(n-1)$, for $S_{n}$ contains an identity submatrix $I_{n-1}$
(we use $S_{n}$ to refer to both the simplex and its matrix) . Since $S_{n}$
has a vector with $-(n-3)$ at the last coordinate, $S_{n}^{0}$ does not have
vectors with the absolute value of the last coordinate different from $0$ or 
$1$; meanwhile, the last coordinate cannot be negative, since it would imply
that one of the first $(n-1)$ coordinates is negative. Thus, $0$ and $1$ are
the only choices for the last coordinate of a vector of $S_{n}^{0}$. If a
vector of $S_{n}^{0}$ has $0$ at the last position, it can have at most one $%
1$ among the other coordinates. Evidently, if a vector of $S_{n}^{0}$ has $1$
at the last position, it can have either $(n-2)$ or $(n-1)$ ones among the
other coordinates. Therefore, the dual of $S_{n}$ consists of the following $%
(0,1)$-vectors: (1) all vectors with $0$ at the last position and one $1$
among the first $(n-1)$ coordinates, i.e. $[1,0^{n-2};0]^{n-1}$; (2) all
vectors with $1$ at the last position and $(n-2)$ $1^{\prime }s$ among the
first $(n-1)$ coordinates, i.e. $[1^{n-2},0;1]^{n-1}$; (3) all vectors with $%
1$ at the last position and $(n-3)$ $1^{\prime }s$ among the first $(n-1)$
coordinates, i.e. $[1^{n-3},0^{2};1]^{\binom{n-1}{2}};$ (4) the zero vector $%
[0^{n}]$. It is easy to see that $(S_{n}^{0})^{0}=S_{n}\cup \lbrack
0^{n-1};1]$, i.e., the double dual of $S_{n}$ is obtained by adding a vector
with zero coordinates, except for only $1$ at the very last position. Denote
by $D$ the matrix whose columns are the elements of $S_{n}^{0}\backslash
\lbrack 0^{n}].$ \ The images of vectors of $S_{n}^{0}\backslash \lbrack
0^{n}]$ under the Voronoi mapping are linearly independent in $R^{N}$, since
there is an $N\times (N-1)$ matrix $M$ such that $D^{T}M=I_{n-1}$. We omit
details here; however we will provide them in a more consecutive paper.
Thus, rank one forms corresponding to the vectors of $S_{n}^{0}$ span a cone
of co-dimension $1$ in $\frak{P}(n).$ This is the cone of all PDQFs for
which $conv\;(S_{n}^{0})^{0}$ is a Delaunay cell in $\mathbb{Z}^{n}$; It is
interesting that for $n=6$ this cone coincides with the wall between the
domains of $\phi _{1}^{6}$ ($\phi _{1}^{6}\backsim $ $D_{6}$) and $%
E_{6}^{\ast }$. \ 

By the above theorem $S_{n}$ is a Delaunay simplex in $\mathbb{Z}^{n}$ for
some PDQFs.
\end{proof}

\section{Tame Wall and Interpretability}

There are many instances of coincidence of $L$-walls and perfect walls. In
particular, for $n\leq 6$ all perfect walls are also $L$-walls, although
already for $n=4$ there are $L$-walls which are not perfect walls ( see
Ryshkov, Baranovskii 1976).

\begin{definition}
\ We call a wall between two $N$-dimensional $L$-type domains $P$
interpretable if it is also a wall between two perfect domains. Conversely a
perfect wall is called $L$-interpretable if is also a wall between two $L$%
-type domains.
\end{definition}

Let $R_{n}$ be the repartitioning complex obtained by adding $[0^{n-1};1]$
to $S_{n}$, i.e. $R_{n}=S_{n}\cup \lbrack 0^{n-1};1].$ Set $%
R_{n}^{+}:=R_{n}\backslash \lbrack 0^{n}].$ As shown above, the Voronoi
images of the vectors of $S_{n}^{0}$ span a cone of co-dimension 1 in $%
\mathbb{R}^{N}$. It follows from the definition of $(0,1)$-dual system that
the interior of this cone consists of \emph{all} PDQFs that have $R_{n}$
among their Delaunay cells. Voronoi  showed that any such cone must be an $L$%
-wall. We call this wall $TW(n)$. Let us prove that for any $n>4$ $TW(n)$ is
a wall between the second perfect form $D_{n}$ \ and a new perfect form.
This new form exhibits a very interesting geometric behavior in all
dimensions, but this will be the subject of another paper.

\section{Perfect Wall Tamed by Big Simplex}

\begin{theorem}
\label{Perfect Domain Description}For any $n>4$ the cone $TW(n)$ is a common
wall of the perfect domain of type $D_{n}$ and the domain of perfect form $%
TF_{n}$, where $TF_{n}=(a_{ij})_{n}\;$is defined as follows. For even $n:$%
\begin{equation*}
a_{ii}=1\text{ if }1\leq i\leq n-1;\;a_{nn}=\dfrac{1}{2}n^{2}-\dfrac{7}{2}%
n+7;\;a_{ij}=\dfrac{n-4}{2(n-2)}\;\text{for }i\neq j,\;j\neq n
\end{equation*}
\begin{equation*}
a_{in}=-1\dfrac{n^{2}-6n+10}{2(n-2)}\text{ for }i<n
\end{equation*}

For odd $n:$%
\begin{equation*}
a_{ii}=1\ \text{if}\ 1\leq i\leq n-1;\ a_{nn}=\dfrac{n^{3}-8n^{2}+23n-20}{%
2(n-1)};
\end{equation*}
\begin{equation*}
a_{ij}=\frac{n-3}{2(n-1)}\;\text{for }i\neq j,\;j<n;\text{ \ }a_{in}=(-1)%
\frac{n^{2}-5n+8}{2(n-1)}\text{ for }i<n
\end{equation*}
In lower dimensions: $TF_{5}\backsim \phi _{2}^{5}$ (III-d perfect form of
Korkine and Zolotareff), $TF_{6}\thicksim E_{6}^{\star }$.
\end{theorem}

\begin{proof}[Proof]
To prove that $TW(n)$ is a wall of a perfect domain we have to complement
the vectors of $\mathcal{V}(R_{n}^{+})\;$to a set $P$ of at least $\frac{%
n(n+1)}{2}$ primitive integral vectors such that the Voronoi image of these
vectors defines a hyperplane in $\mathbb{R}^{N}.$ We call elements of $%
P\backslash R_{n}^{+}$ \emph{complimentary vectors}. Below we give two ways
to complement $\mathcal{V}(R_{n}^{+})$ to the set Voronoi images of minimal
vectors for a perfect form.

A) \emph{Second Perfect Form. }

The complimentary vectors are of type $[1,-1,0^{n-3};0].$ Thus they lie in
one of the half-spaces defined by $TW(n)$. Let \ $(f_{ij})$\ be a quadratic
form defined by the following formulae: $f_{ii}=1$ for $i<n-1$, $f_{nn}=1+%
\binom{n-2}{2},$ $f_{ij}=\frac{1}{2}$ for $i\neq j$ and $i<n,$ and $f_{in}=-%
\frac{(n-2)}{2}.$\ The minimum of \ $f$ \ is $1$ and it is attained on all
vectors of $R_{n}^{+}$ and all $\binom{n-1}{2}$ vectors of type $%
[1,-1,0^{n-3};0]$. Denote by $CompD_{n}\;$\ the set of all vectors whose
coordinates are obtained by circular permutations of the first $n-1$
positions of $[1,-1,0^{n-3};0],$ except for $[-1;0^{n-3};1;0].$ Let us show
that $(f_{ij})$ is integrally equivalent to $D_{n}$. With respect to the
scalar product defined by $(f_{ij}),$ the following integral vectors forms a
Coxeter diagram for $D_{n}:$ $CompD_{n},[0^{n-1};1;0],[1^{n-2};0;1]$. In
this diagram $[0^{n-3};1;-1;0]$ is the vertex of valence 3, and $%
[0^{n-1};1;0]$, $[1^{n-2};0;1]$ are the vertices of the two leaves of the
diagram which are adjacent to $[0^{n-3},1,-1;0]$. Vectors $[0^{n-1};1;0],$ $%
[1^{n-2};0;1]$ and vectors of $CompD_{n}$ \ obviously form a basis of $%
\mathbb{Z}^{n}$. Therefore $(f_{ij})$ is integrally equivalent to $D_{n}.$
Notice that $|\{R_{n}^{+}\cup \lbrack 1,-1,0^{n-3};0]^{\frac{n-1}{2}%
}|=n(n-1),$ which is half the number of minimal vectors of $D_{n}$.
Therefore $R_{n}^{+}\cup \lbrack 1,-1,0^{n-3};0]^{\frac{n-1}{2}}$ and their
inverses are all of the minimal vectors of $(f_{ij})\backsim D_{n}$.

$B)$ \emph{Generalization of the Third Perfect Form. }

Notice that $TF_{n}$ is $1$ for all $\mathbf{v\in }R_{n}^{+}$. Set $%
w:=\lfloor \frac{n}{2}\rfloor .$ The choice of complimentary vectors depends
on the parity of the dimension.

\begin{enumerate}
\item  When $n$ is even and $n>4$ the complimentary vectors for $TF_{n}$
are: $[(w-2)^{n-1};w-1],$ $[w-1,(w-2)^{n-2};w-1]^{n-1},[(w-1)^{n-1};w].$ The
total number of minimal vectors of is $n(n+3)$

\item  When $n$ is odd the complimentary vector for $TF_{n}$ is: $%
[(w-1)^{n-1};w].$ The total number of minimal vectors is $n(n+1)$
\end{enumerate}

With respect to the standard scalar product $TW(n)$ is defined by the
equation $\mathbf{n}\bullet \mathbf{x=0,}$ where $\mathbf{n}$ is given by
formulae $n_{ii}=0$ for $i<n;$ $n_{in}=\frac{n-2}{2(n-4)}$ for $i<n$; $%
n_{ij}=1$ for $i\neq j,$ $j<n$ $n_{nn}=\frac{n^{3}-9n^{2}+24n-19}{2(n-4)}.$
For any $n>3$ $TF_{n}$ is a unique hyperplane in $\mathbb{R}^{N}$ containing
both $R_{n}^{+}$ and the complementary vectors, because for both even and
odd $n$ the complementary vectors form non-zero scalar products with $%
\mathbf{n}$ and therefore do not belong to the hyperplane containing the
cone $TW(n).$ To prove that $TF_{n}$ is positive definite it is enough to
show that $\det TF_{n}>0$, since all other main minors of $TF_{n}$
correspond to bases formed by vectors of length $1$ with angles $\arccos \
a_{ij}$\emph{\ }between them (here we use the correspondence between bases
and quadratic form). If we can construct a basis whose Gramm matrix is $%
(a_{ij}),$ then $\det TF_{n}>0$. Evidently, there exist $n-1$ vectors $%
\mathbf{v}_{1},...\mathbf{v}_{n-1}$ in $\mathbb{R}^{n}$ so that $\mathbf{v}%
_{i}\cdot \mathbf{v}_{j}=a_{ij}$, for $i,j<n$. Now, the norm of $\mathbf{v}%
_{n}$ must be $a_{nn}.$ If the $\cos $ of the (acute) angle between, say, $%
\mathbf{v}_{1}$ and $\mathbf{d}:=\sum\limits_{j=1}^{n-1}\mathbf{v}_{j}$ $(%
\mathbf{d}$\textbf{\ }\emph{is the diagonal of the parallelogram based on
basis }$\mathbf{v}_{1},...\mathbf{v}_{n-1})$\emph{\ }is greater than $%
|\arccos \ \dfrac{a_{in}}{\sqrt{a_{ii}a_{nn}}}|,$ where $i<n,$ then one can
construct vector $\mathbf{v}_{n}$ of norm $a_{nn}$ such that $\mathbf{v}%
_{i}\cdot \mathbf{v}_{n}=a_{in}$.

For odd $n$ this is equivalent to showing that: 
\begin{equation*}
\dfrac{a_{in}^{2}}{a_{nn}a_{ii}}=\dfrac{2(n^{2}-6n-10)^{2}(n-1)}{%
4(n-2)^{2}(n^{3}-8n^{2}+23n-20)}<\dfrac{(2+\frac{(n-2)(n-3)}{n-1})^{2}}{%
2(n^{2}-5n+8)}=\dfrac{(\mathbf{v}_{1}\cdot \mathbf{d})^{2}}{(\mathbf{v}%
_{1}\cdot \mathbf{v}_{1})(\mathbf{d}\cdot \mathbf{d})}
\end{equation*}

For even $n$ this is equivalent to showing that: 
\begin{equation*}
\dfrac{a_{in}^{2}}{a_{nn}a_{ii}}=\dfrac{(-10+n^{2}-6n)^{2}}{%
2(n^{2}-7n+14)(n-2)^{2}}<\dfrac{2(n-2)}{(n-1)}=\dfrac{(\mathbf{v}_{1}\cdot 
\mathbf{d})^{2}}{(\mathbf{v}_{1}\cdot \mathbf{v}_{1})(\mathbf{d}\cdot 
\mathbf{d})}
\end{equation*}

Using elementary algebra and calculus we have checked that both these
inequalities hold for all $n>4.$

The arithmetic minimum of $TF_{n}$ is 1 and the number of minimal vectors is 
$n(n+3)$ for even $n,$ and $n(n+1)$ for odd$n$. For example, this can be
shown by the method of prejective inequalities orginated by \textit{Korkin
and Zolotareff} (see Anzin (1991)). Unfortunately the proof is too tedious
and we have to leave it out. We will publish the proof in another, more
technical paper.

For all $n$ the number of minimal vectors of $TF_{n}$ is equal to that of $%
h_{n}$ of Anzin (1991). It is not difficult to show (Anzin, private
communication) that for $n=5,6,7$ $TF_{n}$ is equivalent to his form $h_{n}$%
. 
\end{proof}

We refer to $TW(n)$ as a \emph{tame perfect wall} because it admits an
interesting $L$-interpretation described above. For $n=6$ $TF_{6}\backsim
E_{6}^{\ast }$ (proved in ''Voronoi-Dickson Hypothesis...'', this volume)
and $TW(n)$ is one of the three (up to $GL(n,\mathbb{Z})$-equivalence) walls
of the domain of $E_{6}^{\ast }$ (see \textit{Barnes} (1957) ) and $.$ The
other two walls, called $W_{2}(24)$ and $W_{3}(21)$ by Barnes, are \emph{%
wild,} as there is no change of $L$-type at almost all interior points of
these perfect walls. The proof that $W_{2}(24)$ is not interpretable can be
found in the other paper by us from this volume. We plan to publish the
proof that $W_{3}(21)$ is not interpretable later. For $n=7$ $TF_{7}\backsim
\phi _{15}^{7}$ from Stacey's (1973, 1975) list (see Anzin (1991) and 
\textit{Martinet }(1996) for $\phi _{15}^{7}$). After this paper had been
submitted for publication Maxim Anzin noticed that in lower dimensions ($%
n=5,6,7$) our series coincides with Anzin's (1991) series $h_{n}.$ Maxim
informed us that he is about to prove that $TF_{n}\backsim h_{n}$ for all $n.
$

\section{The case of n=6: $L$-partition of $E_{6}$}

In this subsection we discuss Delaunay tilings of lattices lying in a small
neigbourhood of $E_{6}$ in the space of parameters. More specifically, we
look at the $L$-partition of $\frak{P}(n)$ near the ray corresponding of $%
E_{6}$. The Delaunay tiling of lattice $E_{6}$ is formed by congruent copies
of the Gosset polytope ($2_{21}$in Coxeter's notation), which is the convex
hull of a unique two-distance spherical set in $\mathbb{E}^{6}$. We refer to
the Gosset polytope as the $G_{6}$-tope. The $G_{6}$-topes of the Delaunay
tiling of $E_{6}$ fall into two translation classes. The star of a lattice
point is formed by 54 $\medskip G_{6}$-topes, 27 in each translation class.

The $G_{6}$-tope is quite remarkable. It has 27 vertices, 216 edges, 72
regular simplicial facets, and 36 regular cross-polytopal facets (e.g. 
\textit{Coxeter }(1995)). Thus, the vertices of the $G_{6}$-tope form a
spherical two distance set. Polytopes whose vertices form a spherical two
distance sets are interesting combinatorial objects (see \textit{Deza and
Laurent }(1997), Deza, Grishukhin, Laurent (1992)). In the case of $G_{6}$%
-tope the two distance structure is realized so that for each vertex $%
\mathbf{v}$ of the $G_{6}$-tope there is a vector $\mathbf{p}_{v}$ such that
the vertex set of $G_{6}$-tope can be represented as $\mathbf{v}\cup
V_{1}\cup V_{2}$, where $V_{1}=\{\mathbf{u}\in S\,\mid \,(\mathbf{u-v)}%
\bullet \mathbf{p}=1\}$, and $V_{2}=\{\mathbf{u}\in S\,\mid \,\,(\mathbf{u-v)%
}\bullet \mathbf{p}=2\}.$ For a detailed description of geometric and group
theoretic properties of the $\medskip G_{6}$-tope see (\textit{Coxeter}
(1973, 1995)).

Below, we show that for every subset of vertices of a Dealunay cell of $%
E_{6} $, $E_{6}$ can be perturbed so that this subset becomes a Delaunay
cell for the perturbed lattice. In particular, this implies that there are
perturbations of $E_{6}$ having a Delaunay simplex of volume 3, the maximal
relative volume of a Delaunay lattice simplex in $\mathbb{E}^{6}$.

\begin{proposition}
For every convex polytope $D$ whose vertex set is a subset of the vertex set
of the $\medskip G_{6}$-tope there is a perturbation of $E_{6}$ making $D$ a
Delaunay polytope for the perturbed lattice.
\end{proposition}

\begin{proof}[Proof]
Denote by $\phi _{E_{6}}(x)$ an inhomogenious quadratic function whose
quadratic part is $E_{6}$, and such that $\phi _{E_{6}}(x)=0$ is an
ellipsoid circumscribing the $\medskip G$-tope. For $\alpha >0$ consider
quadratic function 
\begin{equation*}
\phi (\mathbf{x})=\phi _{E_{6}}(\mathbf{x})+\alpha \sum_{\mathbf{v}\notin
D}\;(\mathbf{p_{v}\bullet x}-1)(\mathbf{p_{v}\bullet x}-2).
\end{equation*}
When $\alpha $ is sufficiently small the quadratic part of $\phi _{E_{6}}(x)$
is close to $E_{6}$ in the space of parameters. The ellipsoid $\phi
_{E_{6}}=0$ circumscribes $D$, since forms $\alpha \,(\mathbf{p_{v}\bullet x}%
-1)(\mathbf{p_{v}\bullet x}-2)$, $\mathbf{v}\notin D$ guarantee that all
vertices of the $\medskip G$-tope that are not in $D$ lie outside of $\phi
_{E_{6}}=0$.
\end{proof}

We thank Maxim Anzin for interesting discussion and valuable comments on
this paper.

R. Erdahl

Department of Mathematics \& Statistics,

Queen's University,

Kingston, ON, K7L 3N6, Canada

Email: erdahlr@post.queensu.ca

K. Rybnikov

Department of Mathematics, Cornell University,

Ithaca, NY, 14853, USA

Email: kr57@cornell.edu

\end{document}